\documentclass[twoside,twocolumn, 10pt]{article}

\usepackage[sc]{mathpazo} 
\usepackage[T1]{fontenc} 
\usepackage{microtype} 

\usepackage{mathptmx}
\DeclareMathSizes{10}{10}{8}{6}   

\usepackage[english]{babel} 

\usepackage[hmargin=0.75in, vmargin=0.9in, columnsep=16pt]{geometry} 
\usepackage[small,labelfont=bf]{caption} 

\usepackage{enumitem} 
\setlist[itemize]{noitemsep} 

\usepackage{graphicx}

\usepackage{cite}
\newcommand\notype[1]{\unskip}

\usepackage{amsmath}
\usepackage{amssymb}
\usepackage{bbm}
\usepackage{mathrsfs}
\usepackage{array}

\usepackage[caption=false,font=footnotesize]{subfig}

\usepackage{scrextend}
\usepackage{url}
\usepackage{xcolor}

\usepackage{multirow} 
\usepackage{booktabs,colortbl}

\usepackage{tabularx}
\newcolumntype{Y}{>{\raggedleft\arraybackslash}X}

\usepackage{floatrow}
\DeclareFloatFont{foot}{\footnotesize}
\usepackage{siunitx}
\usepackage[title]{appendix}

\usepackage{verbatim} 

\renewcommand\footnoterule{\vspace*{-3pt}%
	\hrule width 2in height 0.4pt
	\vspace*{2.6pt}}

\renewcommand\footnotemark{}

\usepackage{abstract} 

\usepackage{titlesec} 
\titleformat{\section}[block]{\large\scshape\centering}{\thesection.}{1em}{} 
\titleformat{\subsection}[block]{\large}{\thesubsection}{1em}{} 
\titleformat{\subsubsection}[block]{\itshape}{\thesubsubsection}{1em}{} 

\usepackage{titling} 

\usepackage[hyperfootnotes=false]{hyperref} 
\hypersetup{
	colorlinks,
	linkcolor={red!50!black},
	citecolor={blue!50!black},
	urlcolor={blue!80!black}
}

\pretitle{\begin{center}\Huge\bfseries} 
	\posttitle{\end{center}} 
\title{Effects of Load-Based Frequency Regulation on Distribution Network Operation} 
\author{Stephanie C. Ross,
	Gabrielle Vuylsteke,
	and Johanna L. Mathieu
	\thanks{S. C. Ross and J. L. Mathieu are with the Department of Electrical Engineering
		and Computer Science, University of Michigan, Ann Arbor, MI 48109 USA
		(e-mail: sjcrock@umich.edu; jlmath@umich.edu).}
	\thanks{G. Vuylesteke is with the Department of Mechanical Engineering, University of Michigan, Ann Arbor, MI 48109 USA
		(e-mail: gevuylst@umich.edu).}
	\thanks{This work was supported by the National Science Foundation Graduate Research Fellowship under Grant No. DGE 1256260.}%
}

\date{} 
\renewcommand{\maketitlehookd}{%
	\begin{abstract}
		\noindent This paper examines the operation of distribution networks that have large aggregations of thermostatically controlled loads (TCLs) providing secondary frequency regulation to the bulk power system. Specifically, we assess the prevalence of distribution network constraint violations, such as over- or under-voltages and overloading of transformers. Our goal is to determine the set of constraints that are at increased risk of being violated when TCLs provide regulation. We compare network operation in two cases: first with TCLs operating freely, and second with TCLs controlled to track a regulation signal. Using GridLAB-D, we run power flow simulations of five real distribution networks. Our results indicate that voltage limits are at increased risk of violation when TCLs provide regulation because of increased voltage variation. Effects on transformer aging are more nuanced and depend on the method used for dispatching TCLs. We find that in many distribution networks it may only be necessary to consider voltage constraints when designing a TCL control strategy that protects the distribution network.
	\end{abstract}
}

\begin{document}
	
	\maketitle

\section{Introduction}
Flexible energy resources are the necessary counterpart to renewable generation. When wind or solar power fluctuates, flexible resources are able to adjust their power to maintain the grid's frequency stability. Currently, the types of resources providing flexibility to the grid are shifting; for example, natural gas units and batteries are replacing coal in PJM's regulation market \cite{pjm_pjm_2017}. Aggregated thermostatically controlled loads (TCLs) are a promising flexible resource \cite{callaway_tapping_2009} capable of providing balancing reserves such as frequency regulation, also known as secondary frequency control.

Utilizing loads for regulation requires local distribution networks to transmit the service to the regional power system. As a result, distribution power flows will change, and load actions could cause local constraints, such as thermal limits of components or voltage limits of nodes, to become violated. This is not just a problem with aggregated loads; other distributed energy resources (DERs), such as photovoltaic generation and electric vehicles, can stress distribution networks. In the future, distribution operators will require better sensing and control over networks with high proportions of DERs. A few states in the U.S. have begun initiatives in this vein \cite{pnnl_electricity_2016}: in New York, Reforming the Energy Vision (REV) aims to transform utilities into ``distributed system platform providers'' that would function much like regional transmission organizations but at the distribution level \cite{nypuc_order_2015}; and in Hawaii, regulators are implementing rules that better align DER adoption with grid reliability \cite{hpuc_instituting_2014}. To ensure reliability, it will be essential for a distribution operator to have information on, and some control over, aggregators' actions. For example, an aggregator might be required to notify the operator when planning to participate in a regional market. If the aggregator's actions are likely to cause voltage violations in one area, the operator might restrict the aggregator to DERs in non-vulnerable areas, or the operator might use other DERs to provide voltage compensation. 

\begin{figure*} [t!]
	\centering
	\includegraphics[width=0.8\textwidth]{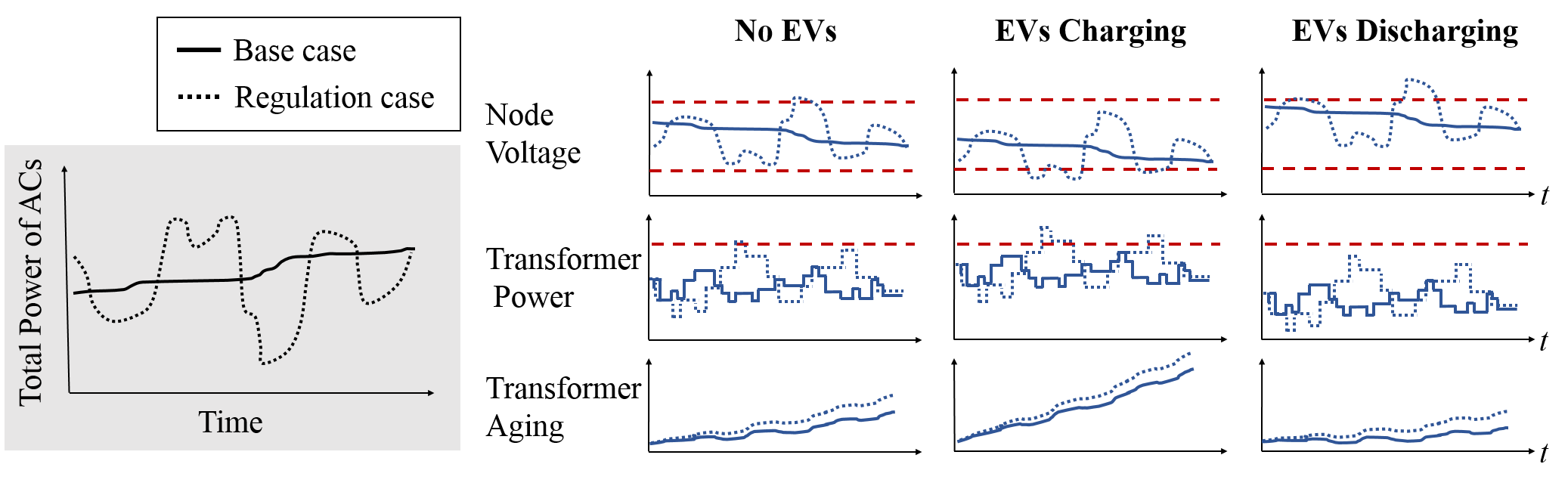}
	\caption{Illustration of possible effects of aggregate DER actions on distribution network operation. In the ``base case'' ACs operate normally, and in the ``regulation case'' ACs are controlled to track a regulation signal. Regulation could cause over/under voltages, as well as high power flows and aging rates in transformers. The addition of EVs to the network could exacerbate constraint violations. Operational limits are indicated with red dashed lines. }
	\label{fig:illustration}
\end{figure*}

In this paper, we assess the impacts of aggregate load control on distribution network operation with the goal of informing coordination strategies between operators and aggregators. Similar studies on network impacts have been conducted for DERs, but typically without control (e.g., residential photovoltaic systems \cite{cohen_effects_2016}, and electric vehicles (EVs) \cite{taylor_evaluations_2010,clement_impact_2010}). Despite a growing literature on regulation provision by TCLs, very little work has considered the impacts of this control on distribution networks. Instead, research has focused on the technical challenges associated with controlling large numbers of distributed resources \cite{callaway_achieving_2011}, \cite{mathieu_state_2013}. To the best of our knowledge, only a few papers have considered the effects of load-based regulation on distribution network operation. The authors of \cite{zhong_relation_2015} found that regulation-provision by TCLs affects distribution transformers' maximum temperatures. The authors of \cite{vrettos_combined_2013} proposed an optimization based approach to provide regulation with TCLs while ensuring distribution constraints are satisfied; this approach is computationally intensive as it requires solving an AC-OPF in each time step. 

As a first step towards a computationally efficient control strategy, this paper identifies the reduced-set of network constraints that are at increased risk of violation when TCLs provide regulation. We identify this set by simulating distribution networks in ``worst case'' load-control scenarios; in these scenarios, TCLs provide their maximum regulation capacity during the network's peak-load hour. The full set of constraints that we assess includes line current, transformer power, transformer aging rate, voltage magnitude, and voltage unbalance.

The contributions of this paper are four-fold: 1) we comprehensively study the effects of load-based regulation on distribution network operation; 2) we examine the effects of having two uncoordinated DER populations (i.e., TCLs and EVs) active on the same network; 3) we explain why some transformers are more likely than others to age faster due to load-based regulation, given our chosen control strategy; and 4) we identify a reduced-set of at-risk network constraints and present preliminary ideas on constraint-protecting control techniques. Contributions 1) and 4) are substantial extensions of our preliminary analysis in \cite{ross_effects_2017}; contributions 2) and 3) are new. 

The paper proceeds as follows. In Section II, the methods are presented. Section III provides the framework for the investigation, which is divided into three separate studies. In Section IV, results are presented and discussed. Section V concludes.

\section{Methods}

\subsection{Illustrative Example}

The illustration in Fig.~\ref{fig:illustration} demonstrates how the actions of aggregated DERs could affect network operation. For instance, voltages may vary more when loads provide regulation, which may increase the prevalence of voltage violations. If EVs are also on the network, their actions (charging or discharging) may improve some constraint violations and worsen others.

\subsection{Simulation Software and Feeder Models}

We use GridLAB-D \cite{gridlabd} to run power flow simulations of distribution networks with time-varying loads. GridLAB-D performs quasi-steady state analysis using a Newton-Raphson algorithm at each time step to solve for a network's three-phase, unbalanced power flow solution. We also use GridLAB-D's dynamic, physics-based models of heating, ventilation, and air conditioning (HVAC) systems. We test networks using a 1 hour simulation with 2 second time steps. To ensure appropriate initialization of dynamic states, we simulate the 24 hours that precede the test hour (using a 30 second time step for the first 23.5 hours and a 2 second time step for the last 0.5 hours).

We use network models from the Pacific Northwest National Lab's (PNNL's) prototypical feeder database \cite{gridlabd}. These models are of actual networks and are prototypical in that they represent common network types in the U.S.~\cite{pnnl_taxonomy}. The networks include fuses, voltage-regulators, capacitor banks, and distribution transformers. Each network is paired with the typical meteorological year (TMY) hourly weather data for its region, which provides the HVAC models with realistic ambient temperature inputs. 

The original feeder models have only one planning load (ZIP model) per distribution transformer; we require higher resolution load models. We disaggregate the transformer-level loads using a method provided by PNNL \cite{pop_script} that estimates the number of houses represented by each planning load and then constructs a model for each house. The house model comprises individual load models (HVAC, water heaters, and pool pumps) and ZIP models that aggregate all other loads in the house. Other modifications made to improve realism are described in Appendix \ref{ap:mods}. 

\subsection{HVAC Modeling}

GridLAB-D's HVAC model includes a thermal model of a house, as well as a model of the space-conditioning devices. The disaggregation method \cite{pop_script} selects randomized parameter values (e.g., wall insulation, house footprint size, temperature setpoints), such that each HVAC model is unique. Throughout this paper we use the term ``AC model'' to refer to the thermal model of the house, as well as the air conditioner device itself. GridLAB-D's AC model is based on \cite{sonderegger_1978} and has three dynamic states: indoor air temperature $\theta$, mass temperature $\theta_\text{m}$, and on/off status $S$ \cite{gridlabd_residential}. These variables evolve according to
\begin{equation} \label{eq:tempDynamics}
\begin{aligned}
\frac{d\theta}{dt} &= \frac{1}{C_\text{a}}\big(-SQ_\text{ac} + rQ_\text{g}-U_\text{a}(\theta-\theta_\text{amb})-H_\text{m}(\theta-\theta_\text{m})\big) \\ 
\frac{d\theta_\text{m}}{dt} &= \frac{1}{C_\text{m}}\big((1-r)Q_\text{g}-H_\text{m}(\theta_\text{m}-\theta)\big) \\
S(t^{\scriptscriptstyle+}) &= \begin{cases} 0, & \theta(t)<\theta_\text{low}\\
1, &  \theta(t)>\theta_\text{high}\\
S(t), & \textrm{otherwise},
\end{cases} 
\end{aligned}
\end{equation}
where $\theta_\text{amb}$,  $Q_\text{g}$, and $Q_\text{ac}$ are time varying inputs that represent outdoor air temperature, internal heat gain, and cooling capacity of the air conditioner, respectively. The constant parameters are defined as: $C_\text{a}$ the thermal mass of indoor air, $C_\text{m}$ the thermal mass of indoor mass, $U_\text{a}$ the overall heat transfer coefficient, $H_\text{m}$ the coefficient of heat transfer between indoor mass and indoor air, $r$ the ratio of the internal gains absorbed by air to the total internal gains, and $\theta_\text{low}$ and  $\theta_\text{high}$ the lower and upper limits of the AC's temperature deadband.

\subsection{External Control of HVACs}

We non-disruptively \cite{callaway_achieving_2011} control ACs with on/off commands that maintain indoor temperature within the user-set deadband. We use the simplest form of the probabilistic dispatch method (see \cite{zhang_aggregate_2013,kara_moving_2013,mathieu_state_2013}) in which a switching probability is broadcast to all ACs and, based on this probability, each AC individually decides whether to switch. Specifically, when the switching probability $u$ is broadcast at the $k$th time step, the $i$th AC will switch if it is available for external control and $p_i(k) < u(k)$, where $p_i(k)$ is drawn from the uniform distribution $\mathcal{U}(0,1)$. An AC is ``available'' for external control if three conditions are satisfied: 1) its indoor air temperature is within the temperature deadband; 2) the AC has not switched within the last two minutes; and 3) if the AC were switched, it would not reach its temperature limit in under two minutes, as predicted by a model. The first condition ensures the control will be non-disruptive to the user and the last two conditions protect the unit's compressor by preventing excessive switching.

We use the simple proportional control scheme from \cite{mathieu_state_2013} to calculate the switching probability $u$. For a desired power level $P_\text{des}$, the switching probability is calculated as 
\begin{equation}
\label{eq:Pcontrol}
u(k+1) = K_\text{P}\frac{(P_\text{des}(k+1)-P_\text{meas}(k))}{\bar{P}_\text{ON}N_\text{AC}(k)}, 
\end{equation}
where the parameter $K_\text{P}$ is a proportional gain, $P_\text{meas}$ is the measured power of the AC population, $\bar{P}_\text{ON}$ is the average power of ACs that are on in steady-state, and $N_\text{AC}$ is the number of ACs that are available to switch. We assume that $N_\text{AC}$ and $P_\text{meas}$ are measured perfectly.

\subsection{Aging Model of Distribution Transformers}

\begin{figure}
	\centering
	\includegraphics[width=\columnwidth]{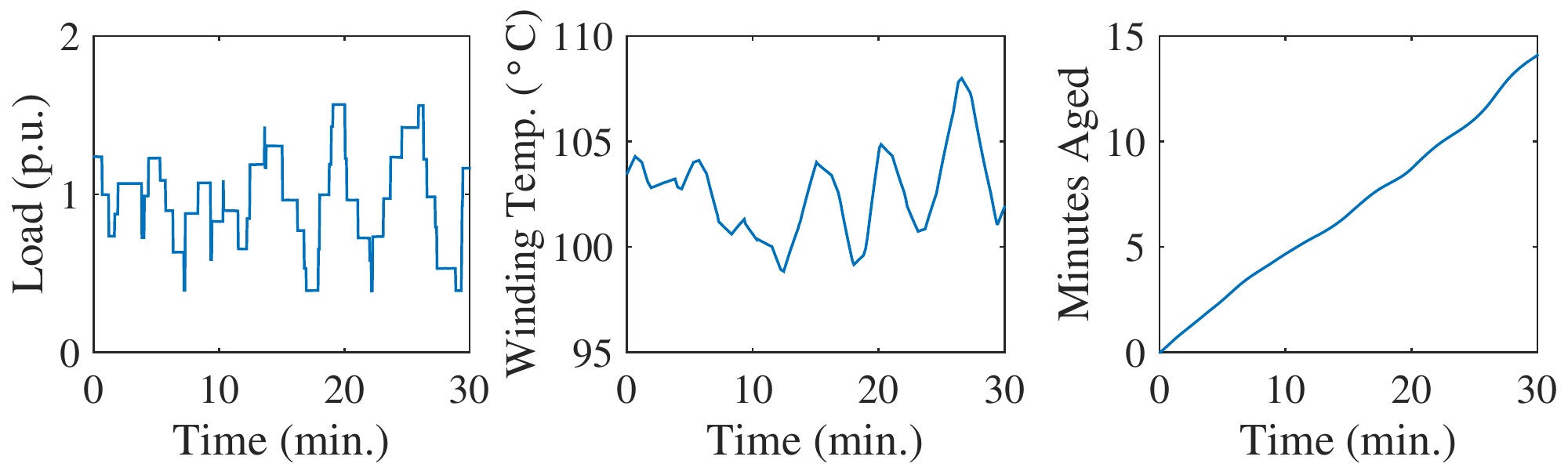}
	\caption{Variables within the transformer aging model: (left) load served by transformer; (middle) temperature of the winding's hot-spot; (right) minutes aged by the transformer. Note that minutes aged are less than the number of simulated minutes because the winding temperature is always $<$110$^\circ$C. }
	\label{fig:xfmr}
\end{figure}

The primary aging mechanism for distribution transformers is the deterioration of coil insulation due to heat from resistive losses \cite{mineraloil}. To estimate transformer aging, we use GridLAB-D's built-in model, which is based on sections 5-7 of IEEE Standard C57.91-1995 \cite{mineraloil}. Figure~\ref{fig:xfmr} shows three key variables within the model: the transformer's load, winding temperature, and estimated minutes aged. 

The model has two dynamic states: $\tilde\theta_\text{oil}$ the difference between the transformer's top-oil temperature and ambient temperature (i.e., $\theta_\text{oil}-\theta_\text{amb}$); and $\tilde\theta_\text{w}$ the difference between the hot-spot winding temperature and the top-oil temperature (i.e., $\theta_\text{w}-\theta_\text{oil}$). Given these states, a transformer's thermal dynamics are described by
\begin{equation}
\begin{aligned}\label{eq:xfmrTherm}
\frac{d\tilde\theta_\text{oil}}{dt} &= \frac{1}{\tau_\text{oil}} (\tilde\theta_\text{oil,u} - \tilde\theta_\text{oil}) \\ 
\frac{d\tilde\theta_\text{w}}{dt} &= \frac{1}{\tau_\text{w}} (\tilde\theta_\text{w,u} - \tilde\theta_\text{w}),
\end{aligned}
\end{equation}
where $\tilde\theta_\text{oil,u}$ and $\tilde\theta_\text{w,u}$ represent the ultimate temperatures that would be reached if the present load were sustained indefinitely \cite{mineraloil}. The ultimate temperatures are computed as  $\tilde\theta_\text{w,u} = \tilde\theta_\text{w,R}L^{1.6}$ and $\tilde\theta_\text{oil,u} = \tilde\theta_\text{oil,R}\big((L^2l_\text{fl}/l_\text{nl} + 1)/(l_\text{fl}/l_\text{nl} + 1)\big)^{0.8}$, where $L$ is the time-varying load (per unit). The oil time constant $\tau_\text{oil}$ is computed according to equations (14) and (15) in \cite{mineraloil}. Values for the thermal parameters in the above equations are derived from data found in \cite{mineraloil,ergon_energy,pacificpowerassoc_2010,vantran_2010} (see Table \ref{tab:xfmrParams}).

\begin{table}
	\caption{Single-Phase Transformer Thermal Parameters}
	\label{tab:xfmrParams}
	\noindent
	\centering
	\begin{minipage}{\columnwidth}
		\renewcommand\footnoterule{\vspace*{-7pt}} 
		\begin{tabular*}{\columnwidth}{@{\extracolsep{\fill}} l r l l}
			\toprule
			\textbf{Parameter}          & \textbf{Value\footnote{Values expressed as a range are parameters that depend on the rating of the transformer (5 kVA-175 kVA).}}          & \textbf{Unit} & \textbf{Source}  \\
			\midrule
			Winding time-constant ($\tau_\text{w}$) & 5         &   min              & \cite{mineraloil} \\
			Winding hot-spot rise ($\tilde\theta_\text{w,R}$) & 80        &   \si{\celsius}           &\cite{mineraloil} \\
			Top-oil rise   ($\tilde\theta_\text{oil,R}$)    & 60       &  \si{\celsius}              & \cite{ergon_energy}\\
			Full-load loss ($l_\text{fl}$)    &  0.0232-0.0112  & {per unit}       & \cite{pacificpowerassoc_2010}\\
			No-load loss  ($l_\text{nl}$)     &  0.0065-0.0042 & {per unit}       & \cite{pacificpowerassoc_2010}\\
			Oil volume         &  5.7-62.7    & {gal}            & \cite{vantran_2010} \\
			Core-plus-coil weight &  56.6-484.9   & {lb}             & \cite{vantran_2010}\\
			Tank-fittings weight  &  67.9-581.9   & {lb}             & \cite{vantran_2010} \\
			\bottomrule
		\end{tabular*}
	\end{minipage}
\end{table}

The final stage of the model calculates a transformer's aging rate $F_\text{AA}$ according to the empirically derived formula \cite{mineraloil}  
\begin{equation} 
\label{eq:ansi}
F_\text{AA} = \exp\bigg(\frac{1500}{383}-\frac{1500}{\theta_\text{w}+273}\bigg),
\end{equation}
where $\theta_\text{w}$ is the winding temperature and $\theta_\text{w} = \theta_\text{amb}+\tilde{\theta}_\text{oil}+\tilde\theta_\text{w}$. The nominal aging rate is equal to one and occurs when $\theta_\text{w} = 110^\circ$C. Note that, in this model, the aging rate is not dependent on the transformer's age. 

\begin{table}
	\caption{Network Constraints}
	\label{tab:constraints}
	\noindent
	\centering
	\begin{minipage}{\columnwidth}
		\renewcommand\footnoterule{\vspace*{-7pt}} 
		\begin{tabular*}{\columnwidth}{@{\extracolsep{\fill}} l l l l }
			\toprule
     		\textbf{Component} & \textbf{Variable} &   \textbf{Lower limit}  & \textbf{Upper limit}   \\
			\midrule
			Service node &  Continuous voltage  &  0.95 p.u. & 1.05 p.u.  \\ 
			&			Emergency voltage    & 0.9 p.u. & 1.083 p.u.   \\
			\noalign{\vskip 0.7mm}   
			3-Phase node & Voltage unbalance &  & 3\%  \\
			\noalign{\vskip 0.7mm}  
			Transformer & Apparent power  &  & 200\% of rating   \\ 
			&	Average aging rate       &   & 1x nominal \\
			\noalign{\vskip 0.7mm}  
		Line & Current  &  & 100\% of rating   \\
			\bottomrule
		\end{tabular*}
	\end{minipage}
\end{table}

\subsection{Assessing the Effects of Regulation}
We test feeders with a scenario designed to reveal the negative effects of load-based regulation. In this scenario, we assume 100\% of residential ACs are controllable, we test networks during the peak-load hour of the year, and we maximize the amplitude of the regulation signal. We find each feeder's peak-load hour by running a simulation over the summer months without regulation. We compare this simulated peak-load to the feeder's original planning load. If the simulated peak is less than 90\% of the planning load, we use \cite{pop_script} to iteratively adjust the simulated peak by repopulating the feeder with more houses until the peak is between 90\% and 100\% of the planning load. For the regulation signal, we select a segment of the PJM Reg-D signal \cite{regd} in which the energy consumed during the regulation hour is equal to that of the base case. We scale the signal such that its maximum magnitude is 0.4 p.u. -- the largest magnitude the AC population can track before performance begins to deteriorate.

We assess the effects of regulation by identifying changes in network variables between the base and regulation cases. Variables of interest and their corresponding constraints are listed in Table \ref{tab:constraints}. First in the table are constraints on the voltage magnitude at service nodes, which are where service lines connect to distribution transformers. The continuous voltage limits are based off of ANSI standard C84.1 \cite{ansi_C84.1} and are only violated if surpassed for over two minutes. The emergency limits are from \cite{gonen_electric_2015} and are similar to those in \cite{ansi_C84.1} and \cite{sce_networklimits}. Second in the table is voltage unbalance of 3-phase nodes: unbalance should be kept below 3\% to keep 3-phase motors from overheating \cite{ansi_C84.1}. 

Third in the table are constraints to prevent transformers from aging too rapidly. Many utilities use apparent power as a proxy for aging rate. We use a limit of 200\% of a transformer's rating as in \cite{sce_networklimits}. Unlike utilities, we have access to our simulated transformers' aging rates. We constrain a transformer's aging rate such that its average rate over the simulation hour must be less than one. If a transformer surpasses this limit, its lifetime will be shorter than nominal (20 years). Fourth in the table is current flow, which is monitored for overhead and underground lines on the primary side of distribution transformers. We set the over-current constraint to 100\% of the line's rating. We also monitor the status of all fuses.

\section{Investigation Framework}

\begin{table}
	\caption{Feeder Characteristics and Peak-Hour Average Conditions}
	\label{tab:feeder}
	\noindent
	\centering
	\begin{minipage}{\columnwidth}
		\renewcommand\footnoterule{\vspace*{-7pt}} 
		\begin{center}
			\begin{tabular*}{\columnwidth}{@{\extracolsep{\fill}} l l c c c c} 
				\toprule
				\multirow{2}{1.3cm}{\textbf{Feeder}} & \multirow{2}{1.5cm}{\textbf{Areas served}} &   \multirow{2}{1cm}{\textbf{Houses with AC}} &
				\multirow{2}{0.7cm}{\centering \textbf{Temp.} (\si{\celsius})} & 
				\multirow{2}{0.8cm}{\centering \textbf{Load} (MVA)} & \multirow{2}{0.65cm}{\textbf{AC Load}}  \\ \\
				\midrule
				R1-12.47-1    & Suburb, rural    & 43\% & 34.7 & 5.71 & 43\% \\
				R2-12.47-2       & Suburb        & 75\% & 34.4 &  7.08  & 51\%  \\
				R3-12.47-3     & Suburb          & 50\% & 45.6 &  8.13 & 42\%  \\
				R4-12.47-1      & Urban, rural   & 95\% & 35.8 & 6.12 & 51\%  \\
				R5-25.00-1    & Suburb, urban   & 97\%  & 33.1 & 8.25 &  55\% \\
				\bottomrule
			\end{tabular*}
		\end{center}
	\end{minipage}
\end{table}

\subsection{Main Study}
We conduct a survey across five networks with the goal of determining how sensitive the results are to different feeder parameters and topologies. Each feeder is from one of five climate regions in the U.S. and varies in terms of voltage level, topology, geographical density of buildings, etc. \cite{pnnl_taxonomy}. Table \ref{tab:feeder} lists the feeders' key characteristics and conditions during the peak-hour. The first two parts of a feeder's name indicates (climate region number)-(voltage level). For brevity, we will refer to the feeders only by region number. In this study, we expect the effects of regulation to be larger for feeders with higher percentages of AC load (i.e., R2, R4, and R5). We also anticipate that regulation could cause voltage issues on feeders with long lines, which are typically in rural areas (i.e., R1 and R4).  

\subsection{Electric Vehicles Study}
In this study, we examine the effects of load-based regulation when EVs are also active on the network. We assume 20\% of residential houses have EVs. All EVs are identical and are able to both charge and provide power to the grid by discharging. We consider two scenarios: ``EV+'' in which all EVs charge for an hour, and ``EV-'' in which all EVs discharge for an hour, where the latter could occur if the EVs were coordinated by an aggregator. We model EVs as constant power loads that charge/discharge at 3.3kW with unity power factor. We conduct this study on feeder R1.

\subsection{Randomization Study}

This study is designed to investigate one of the results from the main study -- that some transformers have an increased aging rate due to load-based regulation. Our goal is to determine whether some transformers consistently have an increased aging rate across multiple trials. We run six randomized trials, all on feeder R1, with different random instantiations of the initial on/off status of the ACs as well as the ACs' probabilistic responses; all other parameters are kept constant.

\section{Results}

\subsection{Main Investigation}\label{sec:mainResults}

\subsubsection{Service-Node Voltage Results} \label{sec:voltresults}

\begin{figure}
	\centering
	\includegraphics[width=\columnwidth]{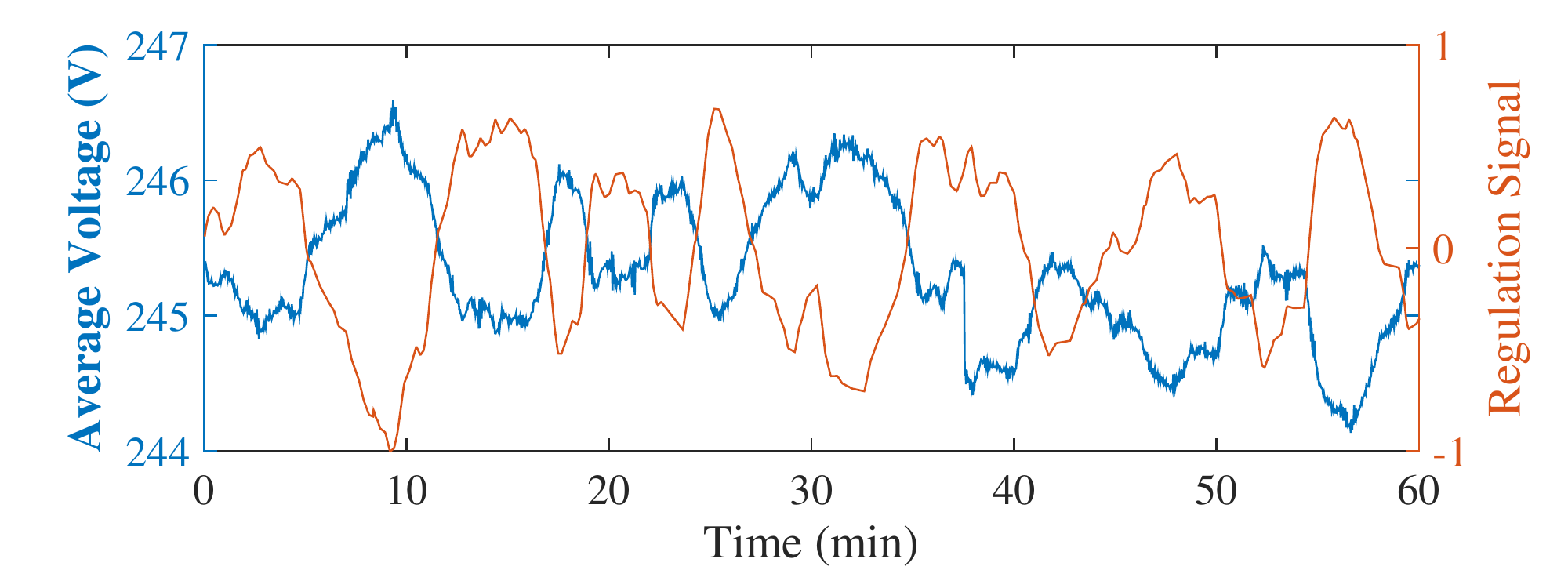}

	\caption{Negative correlation between the regulation signal and average voltage during a regulation case simulation. The voltage trajectory is the average across all residential service nodes. (Data from R1.)}
	\label{fig:voltageCorr}
\end{figure}

\begin{figure}
	\centering
	\includegraphics[scale = 0.3]{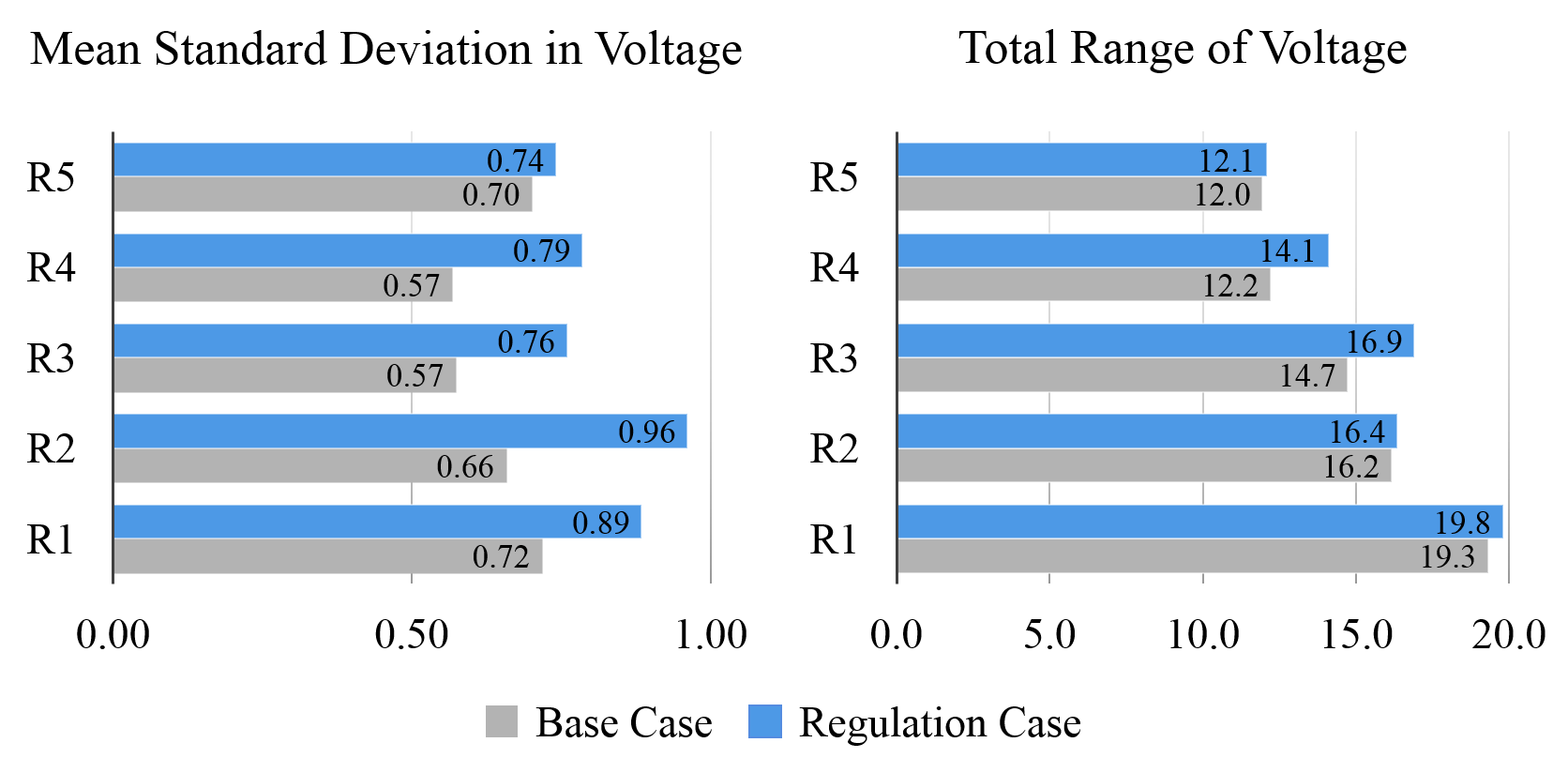}
	\caption{Variation of voltage at residential service nodes. Mean standard deviation of voltages is the mean across all nodes; total range of voltages is the range across all nodes. Both metrics increase in regulation cases for all five networks.}
	\label{fig:voltageBars}
\end{figure}

We find that load-based regulation causes voltages to vary more, but the increase is not large enough to cause constraint violations in the networks we study. As shown in Fig.~\ref{fig:voltageCorr}, voltages deviate in the opposite direction of the change in power injections due to load-based regulation. Figure~\ref{fig:voltageBars} shows that, for all feeders, load-based regulation causes an increase in the variation of voltages at service nodes. Despite the increase in voltage variation, no limit violations occur as a result of regulation for any of these feeders. Feeder R4 has voltages that surpass the continuous upper limit in the regulation case, but for less than the 2 minute duration required for continuous limit violations. For all feeders, there are no violations of emergency limits.

\begin{figure}
	\centering
	\includegraphics[width=\columnwidth]{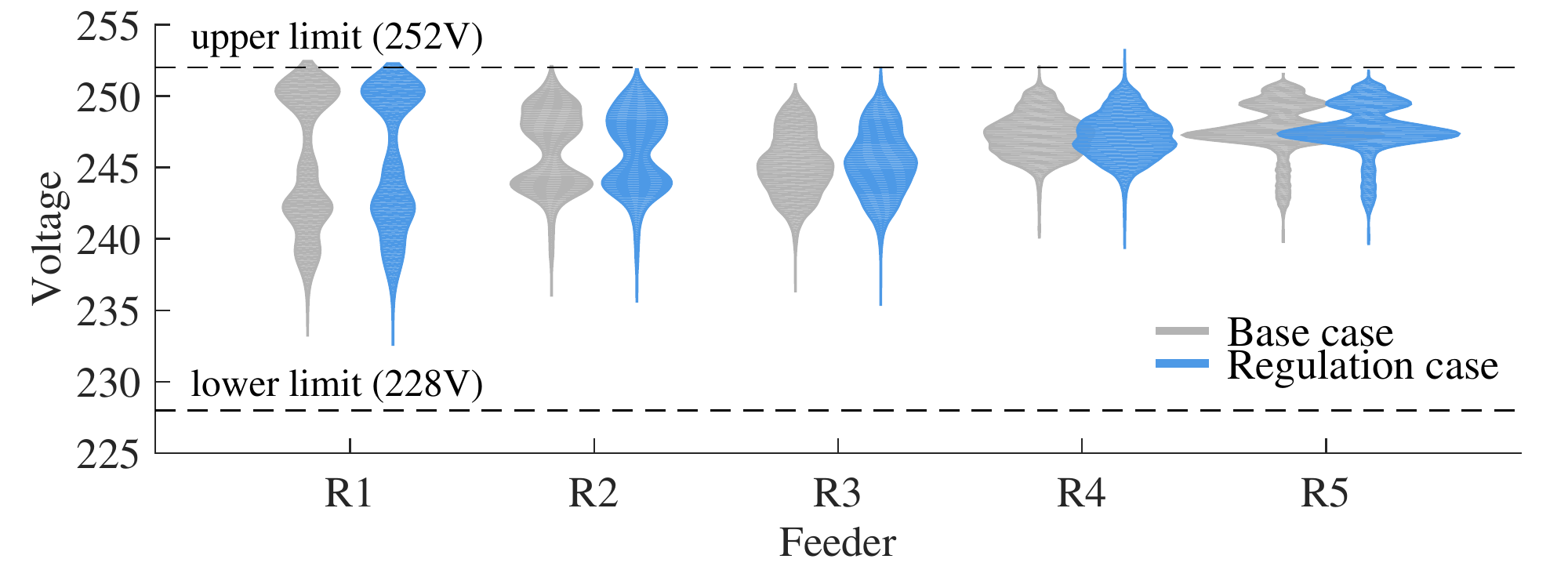}
	\caption{Voltage distributions of residential service nodes. Distributions for R1 and R4 cross the continuous upper limit in both cases. Regulation reduces over-limit voltages for R1 but increases over-limit voltages for R4. Distributions are visualized using kernel density plotting.}
	\label{fig:voltageAllFeeders}
\end{figure}

Figure~\ref{fig:voltageAllFeeders} shows the voltage distributions of all service nodes for the five feeders. An increase in variation due to regulation is indicated by elongation in the distributions. We find that the voltage density outside of (or close to) voltage limits generally increases due to regulation: three of five distributions move closer to the upper limit, and five of five distributions move closer to the lower limit. There would likely be worse impacts on long, rural feeders with poor voltage regulation. In these ``weak'' feeders, we expect voltage distributions would be longer-tailed, and more density would shift outside of the limits in regulation cases. Unfortunately, we did not have access to a weak feeder model for testing.

\subsubsection{3-Phase Node Voltage-Unbalance Results}

Load-based regulation has little impact on voltage unbalance. Although feeders R1 and R4 have nodes that violate the upper limit of 3\%, they do so in both the base case and the regulation case. The largest increase in unbalance due to load-based regulation is only +0.076\%. Its minimal effect on voltage unbalance is likely due to the even distribution of ACs across the three phases of each network.

\subsubsection{Transformer Results}\label{sec:xfmrResults}
 
\begin{table}
	\renewcommand{\arraystretch}{1.1}
	\caption{Transformer Results: Apparent Power and Aging Rate}
	\label{tab:results}
	\noindent
	\centering
	\begin{minipage}{\columnwidth} 
		\renewcommand\footnoterule{\vspace*{-7pt}} 
		\begin{center}
			\begin{tabular*}{\columnwidth}{@{\extracolsep{\fill} } l c c c S[table-format=1.2] S[table-format=1.2]}
				\toprule
				\noalign{\vskip 0.5mm} 
				\multirow{3}{*}{\hspace{-3pt}\raggedleft \textbf{Feeder} } & \multicolumn{3}{c}{\textbf{Population mean}} & \multicolumn{2}{c}{\textbf{\% of pop. w. violation}} \\
				\cmidrule(lr){2-4} \cmidrule(l){5-6} 
				& {Power (p.u.)} & \multicolumn{2}{c}{Aging rate} & {Power} & {Aging} \\ 
				\noalign{\vskip -0.25mm} 
				\cmidrule(lr){2-2}  \cmidrule(lr){3-4} \cmidrule(lr){5-5} \cmidrule(l){6-6}
				&   \emph{base} = \emph{reg.}  & {\centering \emph{base}} & {\centering \emph{reg.}}  & \emph{base} = \emph{reg.} & \emph{base} = \emph{reg.}\\
				\midrule
				R1 &    0.42 & 0.026 & 0.025  & 0.17   & 0   \\
				R2 &   0.55  & 0.076 & 0.070 & 0   & 1.04   \\
				R3 &   0.23  & 0.077 & 0.076 & 0.06  & 1.23   \\
				R4 &  0.43   & 0.015 & 0.015 & 0  & 0.21  \\
				R5 &   0.47  & 0.031 & 0.030 & 0.54  & 0.81  \\				
				\bottomrule
			\end{tabular*}
		\end{center}
	\end{minipage}
\end{table}

\begin{figure*}
	\centering
	\subfloat[][Transformer apparent power distributions]{
		\includegraphics[width = 0.32\textwidth]{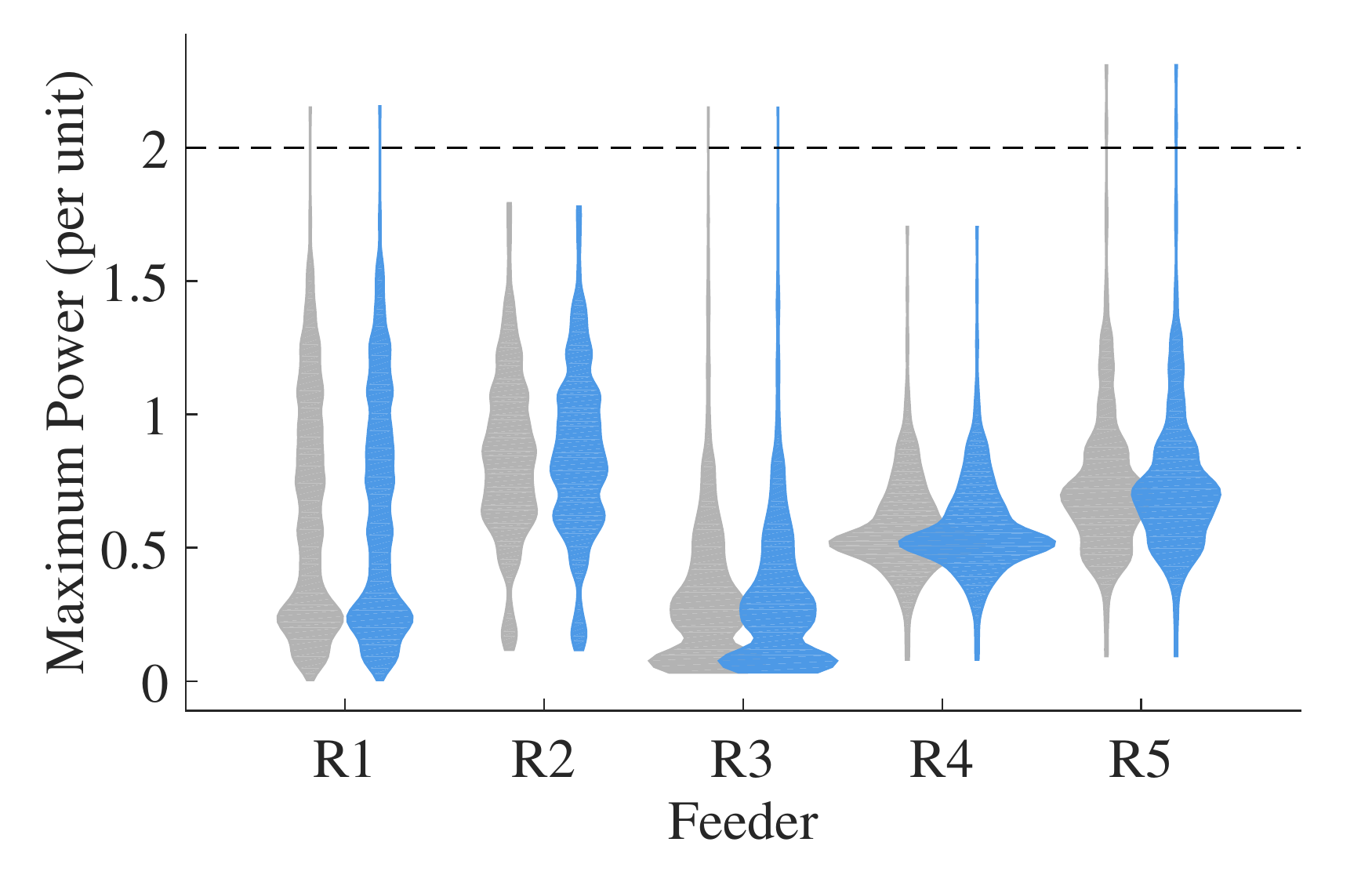}%
		\label{fig:xfmrViolinPwr}
	}
	\subfloat[][Transformer aging rate distributions]{
		\includegraphics[width = 0.32\textwidth]{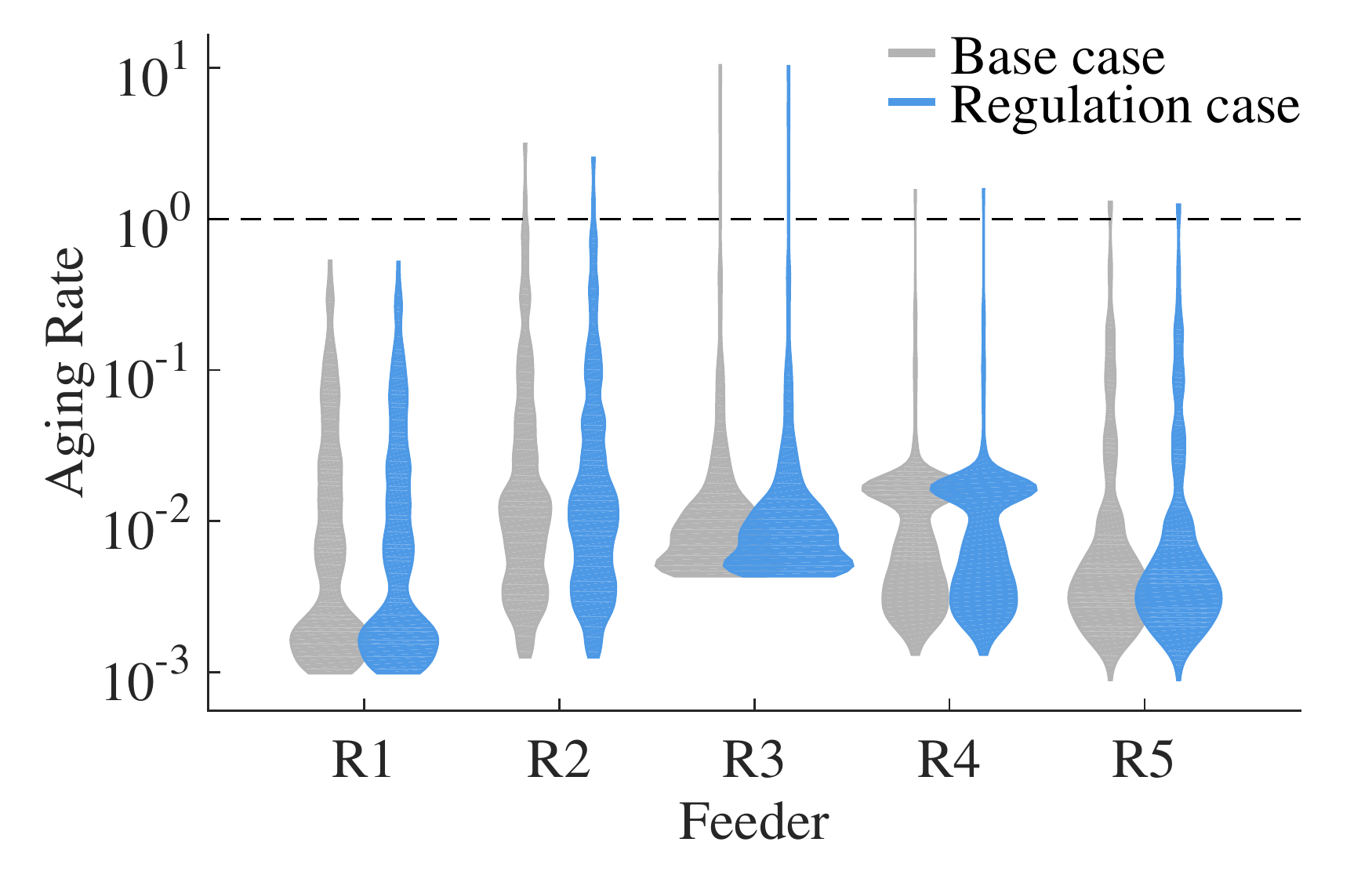}
		\label{fig:xfmrViolinAging}
	}
	\subfloat[][Change in aging rate distributions]{
		\includegraphics[width = 0.32\textwidth]{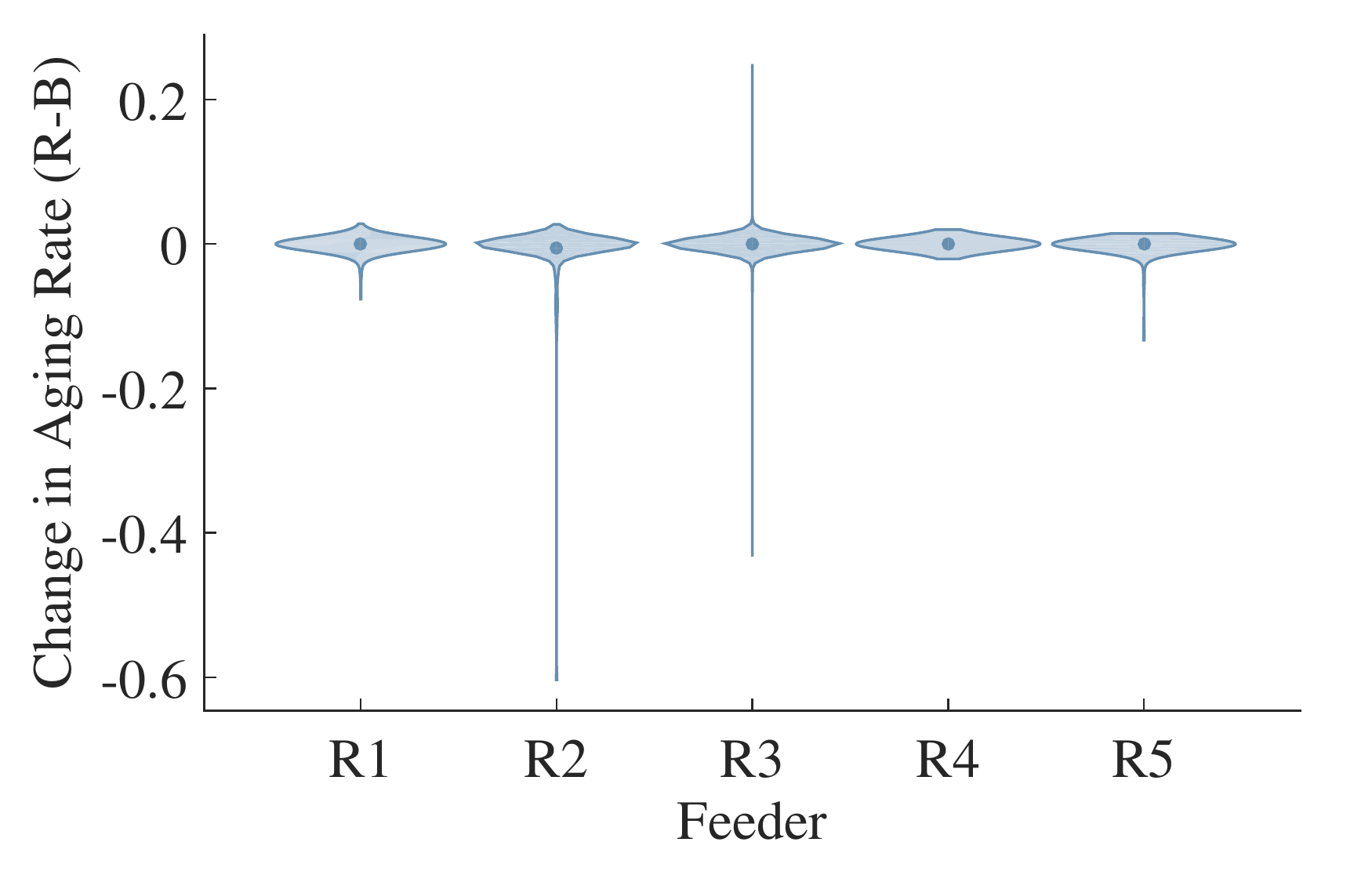}
		\label{fig:xfmrViolinChange}
	}

	\caption{Transformer population distributions for (a) power, (b) aging rate, and (c) change in aging rate (regulation - base). Plots (a) and (b) show that the power and aging rate distributions are positively skewed and do not change substantially from base case to regulation case. Plot (c) shows that the change in aging rate is mostly symmetric about the $y=0$ line, excluding a few large outliers. }%

	\label{fig:xfmrViolin}%
\end{figure*}

For all networks, we find that the transformers experience very few apparent power or aging rate violations and that load-based regulation has little effect on these violations as shown in Table~\ref{tab:results}. For all five feeders, the per-unit apparent power of transformers, averaged across the population, is less than 0.6 p.u. (where the p.u. base is the transformer rating) and does not change from base case to regulation case. These low loading levels translate to low mean aging rates. In four of five feeders, the mean aging rate decreases slightly due to regulation, indicating that regulation may be beneficial for some transformers. The percentage of the population with violations (apparent power or aging rate) is small and unaffected by regulation.

Distributions of the transformer populations' loading and aging rates are shown in Fig.~\ref{fig:xfmrViolin}. In each feeder, a small number of transformers have much higher loading and aging rates than the rest of the population, as indicated by the positive skews in Fig.~\ref{fig:xfmrViolin}\subref{fig:xfmrViolinPwr} and \subref{fig:xfmrViolinAging}. Notably, the aging rate distributions require a log-scale for visualization. Figure~\ref{fig:xfmrViolin}\subref{fig:xfmrViolinChange} demonstrates that the changes in aging rate from base to regulation case are generally small. As indicated by the thin spikes, there are a few transformers that experience large changes in aging rate. Excluding these outliers, the distributions in Fig.~\ref{fig:xfmrViolin}\subref{fig:xfmrViolinChange} are approximately symmetric about the line $y=0$, i.e., about 50\% of transformers age faster due to regulation, and about 50\% age slower.

\subsubsection{Line Results} In all simulations, line-current stays well below the 100\% conductor rating constraint, and all fuses remain closed.

\subsection{Results: Electric Vehicle Study} \label{sec:evresults}

We find that -- when EVs are charging -- voltage violations due to load-based regulation worsen. Voltage results are presented in Table~\ref{tab:EVresult}. We restrict our attention to the 1.05 p.u. continuous limit because it is the only constraint violated. In the base case, the results for the EV+ and No-EV trials are similar: 5.69\% and 5.35\% of nodes exceed the 1.05 p.u. limit for at least one time step, respectively. However, in the regulation case, over-limit nodes increase to 20.07\% in the EV+ trial and decrease to 4.85\% in the No-EV trial. A similar divergence occurs in the EV+ and No-EV trials when examining nodes with continuous voltage violations (i.e., over-limit for more than 2 minutes). In contrast to the other trials, the EV- trial has only 1.34\% of its nodes over-limit for any duration, and this percentage remains almost constant across cases and durations. 

\begin{figure}
	\centering
	\includegraphics[width = \columnwidth]{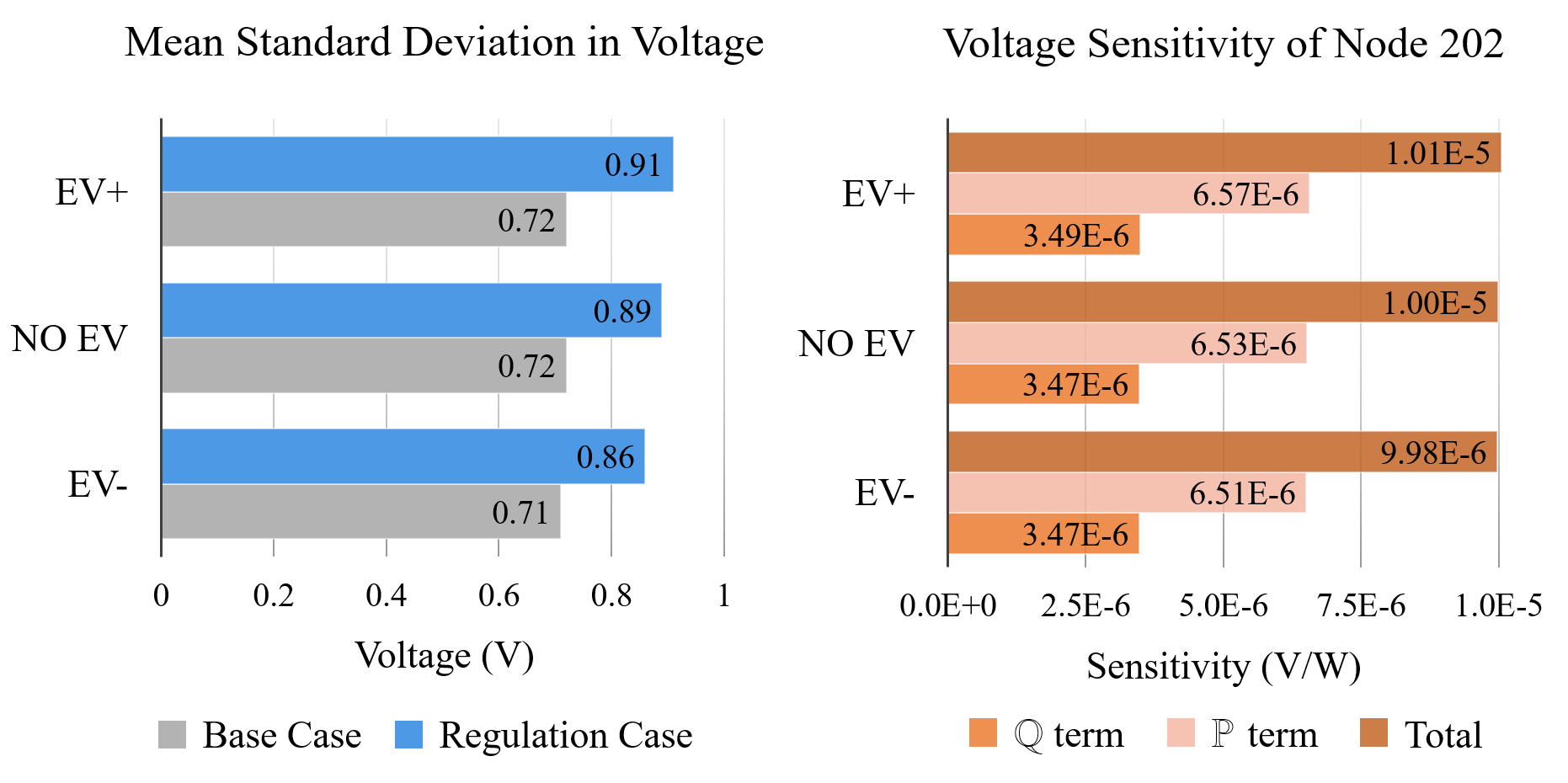}
	\caption{Voltage variation and sensitivity across three EV trials: (left) the change in standard deviation due to regulation increases as loading levels increase across trials (from EV- to EV+) and (right) the computed voltage-sensitivity of a representative node increases as loading levels increase.}
	\label{fig:voltageBarsEV}
\end{figure}

The divergence in the No-EV and EV+ results is due, in part, to network-voltages' increased sensitivity to load variation (due to load-based regulation), at higher levels of base load. Empirically, we see that as base loading increases (from EV- to EV+) the change in voltage variation due to load-based regulation increases (see Fig.~\ref{fig:voltageBarsEV}). Analytically, we find that the sensitivity of voltage to changes in power injections due to regulation increases as loading levels increase. To derive this sensitivity, we use the DistFlow voltage equation \cite{baran_optimal_1989}, as reformulated in \cite{cespedes_new_1990}: 
\begin{equation}\label{eq:distflow}
V_k^4 + (2(rP_{ik}+xQ_{ik})-V_i^2)V_k^2 + (r^2+x^2)(P_{ik}^2+Q_{ik}^2) = 0,
\end{equation}
where $V_i$ and $V_k$ are the voltage magnitudes at buses $i$ and $k$, $P_{ik}$ and $Q_{ik}$ are the real and reactive power flows received by bus $k$, and $r$ and $x$ are the resistance and reactance of the line connecting the buses. Solving the quadratic for $V_k^2$ and then taking the square root, we obtain an expression for $V_k$. Given this expression, the sensitivity of $V_k$ to changes in power injection at bus $k$ due to regulation is
\begin{equation} \label{eq:partials}
\frac{\partial V_k}{\partial P_\text{reg}} = \underbrace{\frac{\partial V_k}{\partial Q_{ik}}\frac{\partial Q_{ik}}{\partial P_\text{reg}}}_{\displaystyle \mathbb{Q}} + \underbrace{\frac{\partial V_k}{\partial P_{ik}}\frac{\partial P_{ik}}{\partial P_\text{reg}}}_{\displaystyle \mathbb{P}}+ \underbrace{\frac{\partial V_k}{\partial V_{i}}\frac{\partial V_{i}}{\partial P_\text{reg}}}_{\displaystyle \mathbb{V}},
\end{equation}
where $\mathbb{Q}$, $\mathbb{P}$, and $\mathbb{V}$ are the three sensitivity terms. To evaluate \eqref{eq:partials}, we must find expressions for each derivative. We assume $V_i$ is an infinite bus; thus $\partial V_i/\partial P_\text{reg} = 0$ and $\mathbb{V} = 0$. We solve for the partials of $V_k$ in the $\mathbb{Q}$ and $\mathbb{P}$ terms using the previously derived expression for $V_k$. To solve for the remaining derivatives, we assume that changes in flow are only due to changes in regulation power, i.e., $\partial P_{ik}/\partial P_\text{reg}=1$ and $\partial Q_{ik}/\partial Q_\text{reg}=1$. Finally, we find $\partial Q_{ik}/\partial P_\text{reg}$ given that ACs have a constant power factor of 0.97.

Using \eqref{eq:partials}, we find an expression for the voltage sensitivity at the downstream node of a representative line in R1. We evaluate the node's voltage sensitivity at all operating points from the base case simulation. For each trial, the computed voltage sensitivity, averaged across time, increases as base loading increases (see Fig.~\ref{fig:voltageBarsEV}). Note that $\mathbb{P}$ contributes more to the increase in sensitivity than $\mathbb{Q}$.

The increase in voltage sensitivity with increased loading may explain why there is a large change in over-limit nodes in the regulation case of the EV+ trial, but only small changes in the other trials (see Table \ref{tab:EVresult}). It should be noted that the EV- trial results deviate from what might be expected (i.e., an increase in voltages) because a capacitor bank switches off due to an initial rise in voltage; the end result is a net decrease in voltage for most nodes. 

\begin{table}
	\caption{EV Trials: \% of Nodes above 1.05 p.u. Voltage Limit}
	\label{tab:EVresult}
	\noindent
	\centering
	\begin{minipage}{\columnwidth}
		\renewcommand\footnoterule{\vspace*{-7pt}} 
		\begin{tabular*}{\columnwidth}{@{\extracolsep{\fill}} l l S S S S}
			\toprule
			\noalign{\vskip 0.25mm} 
			&  & \multicolumn{2}{c}{\textbf{Any duration}}  & \multicolumn{2}{c}{\textbf{More than 2 min}}  \\ 
			\noalign{\vskip -0.25mm} 
			\cmidrule(r){3-4} \cmidrule(l){5-6}
			& & \emph{base} & \emph{regulation} & \emph{base} & \emph{regulation} \\
			\noalign{\vskip -0.25mm} 
			\midrule
			EV+ & charge & 5.69 & 20.07 & 0.17 & 0.84  \\
			No-EV & n/a & 5.35 & 4.85 & 0.17 & 0  \\
			EV- & discharge & 1.34 & 1.34 & 1.34  & 1.17  \\
			\bottomrule
		\end{tabular*}
	\end{minipage}
\end{table}

\subsection{Results: Randomization Study}

\begin{figure}
 	\centering
	\includegraphics[width = 0.8\columnwidth]{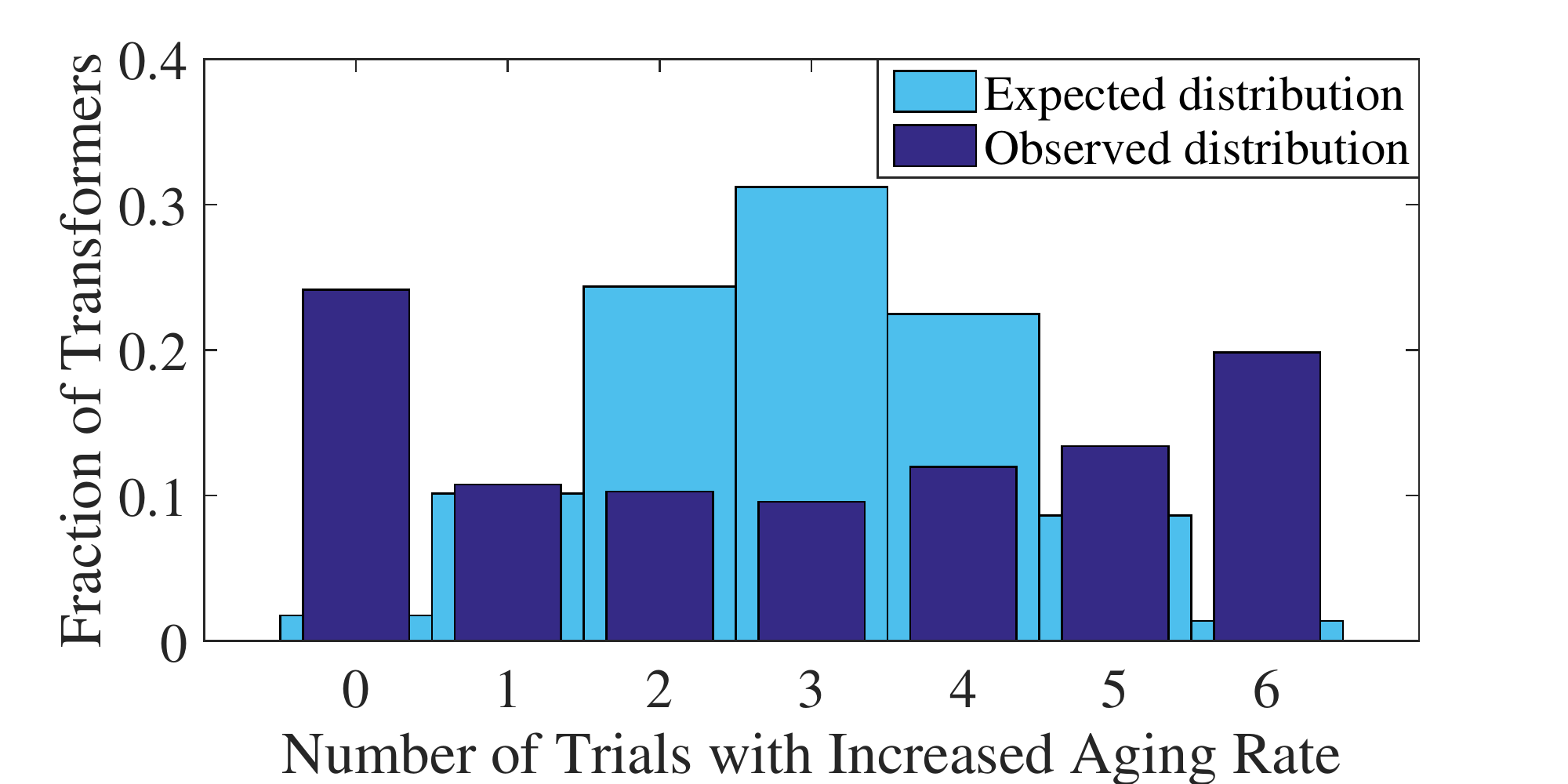}%

 	\caption{Observed and expected distributions for the number of randomized trials a transformer will experience an increased aging rate due to regulation. The difference between the two distributions indicates that some transformers are more likely than others to experience an increased aging rate.}%
 	\label{fig:randAging}%
 \end{figure}

We find that some transformers are more likely to age faster due to load-based regulation and some are more likely to age slower. Figure~\ref{fig:randAging} shows the observed distribution of the number of trials transformers experience an increased aging rate due to regulation. We compare this distribution to the ``expected distribution'', which represents the hypothesis that all transformers are equally likely to have an increased aging rate. Specifically, we model the transformers' aging outcomes as i.i.d. Bernoulli random variables with probability 0.49 of increased aging. (This value is the observed prevalence of increased aging in the randomized trials.) The expected distribution, then, is the number of increased aging outcomes that occur over six Bernoulli trials, i.e., a Binomial distribution. The difference between the two distributions in Fig.~\ref{fig:randAging} indicates that the aging outcomes cannot be modeled as i.i.d. random variables. Instead, the edges of the observed distribution indicate that a large portion of transformers is very likely to experience an increased aging rate, and another large portion is very likely to experience a decreased aging rate. 

\begin{figure}
	\centering
		\includegraphics[width = 0.8\columnwidth]{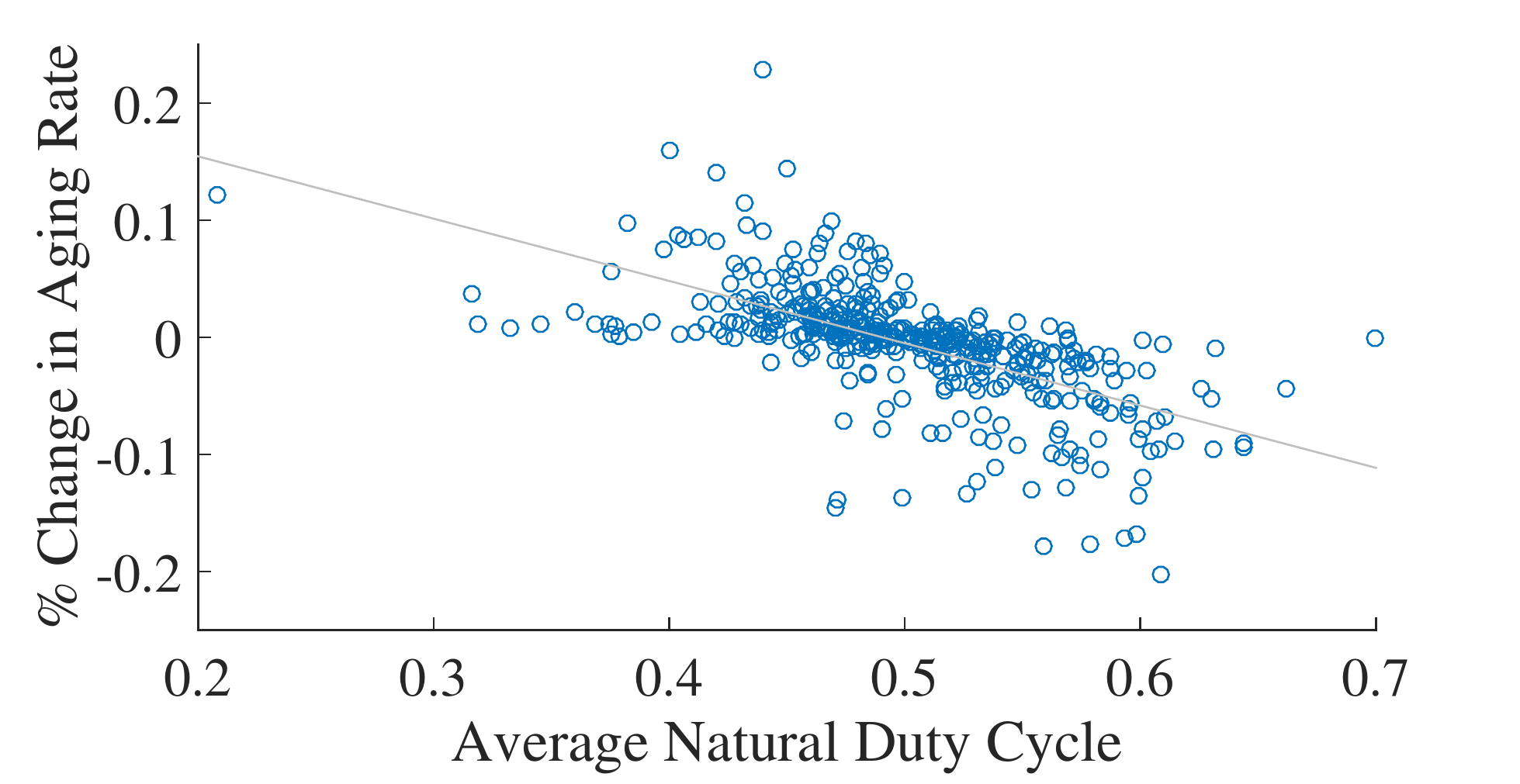}
		\caption{Scatter plot of each transformer's percent change in aging rate (averaged across the randomized trials) versus the average natural duty cycle of the ACs supplied by the transformer. The data is negatively correlated, with correlation coefficient = 0.628.}
		\label{fig:dutyCycleCorr}
\end{figure}

Load-based regulation does not have a consistent effect on transformers' aging rates, in part, because of an unanticipated effect of probabilistic dispatch. In the following discussion, we will refer to switching in response to probabilistic dispatch as ``dispatch-switching''. We will also refer to a cycle through an AC's deadband without dispatch-switching as a ``natural cycle''. When all ACs receive the same $u$, an AC with a low natural duty-cycle is more likely to be dispatch-switched on than off, and an AC with a high natural duty-cycle is more likely to be dispatch-switched off than on. See Appendix~\ref{appendix:duty} for a derivation of this result. If an AC is dispatch-switched on more than off, it will ``cool-cycle'', i.e., cycle in the cooler part of its deadband. This cycling behavior reduces the average temperature of the house, thus increasing its average power consumption and, consequently, the transformer's aging rate. A similar argument can be made for why warm-cycling reduces transformer aging rates. As shown in Fig.~\ref{fig:dutyCycleCorr}, there is a negative correlation between the percent change in aging rate of a transformer and the average natural duty cycle of the ACs the transformer supplies. Transformers supplying ACs with a lower than average duty cycle ($<0.503$) frequently have an increased aging rate, and transformers supplying ACs with a higher than average duty cycle frequently have a decreased aging rate.

\subsection{Control Recommendations}

Our first step towards designing a control strategy for load-based regulation that is both network-protecting and computationally efficient is to identify a reduced-set of network constraints. Although this set will be network dependent, we expect that voltage unbalance and line current constraints could be safely omitted from the reduced-set for most networks. Transformer aging and power constraints could also be omitted, since load-based regulation did not cause any additional violations in our results. However, some utilities may choose to enforce stricter constraints that protect transformers from increased aging rates. Finally, we expect that in most networks we will need to include continuous voltage constraints in the reduced-set because of the increase in voltage variation due to regulation and the increased risk of violations in high loading scenarios.

We propose a few strategies for providing load-based regulation while protecting distribution networks. First, to prevent cool-cycling and thereby the likelihood of increased aging rates for transformers, ACs closer to naturally switching should be dispatched first (i.e., priority-stack control \cite{hao_aggregate_2015}). Second, to prevent voltage violations, the amount of participation by loads in voltage-weak areas could be reduced, or separate resources could be dispatched for voltage support. Reduction in participation could be implemented in one of two ways: 1) by allowing only a portion of loads in the area to participate in regulation, or 2) by dispatching loads in that area less frequently. Both of these strategies require some coordination between DER aggregators and the distribution operator.

\section{Conclusions}

In this paper, we identified the set of distribution network constraints at risk of increased violation due to the provision of regulation by TCLs. Ultimately, this reduced-set of constraints is network-dependent; however, in the five networks studied here, load-based regulation had fairly benign effects on distribution network operation, which suggests that the constraint-set could be reduced substantially. Based on our results, we have recommended excluding line current and voltage unbalance constraints, possibly including transformer aging constraints, and almost always including voltage constraints. 

It is possible that the impacts of regulation were not severe in the networks studied here because of their typical nature. In the future, we plan to investigate the impacts of regulation on non-typical feeders, particularly ones with voltage-weak areas. Additionally, future distribution networks may be much more active and complex than those modeled here. For instance, a network could simultaneously have high penetrations of photovoltaics, EVs, and regulation-providing loads, as well as multiple aggregators controlling portions of each DER population. Our future work will address the distribution operator's role in these types of networks -- in particular, how to equitably regulate DER aggregators' actions while ensuring the reliability of the distribution network.

\begin{appendices}

\section{Feeder Model Modifications}\label{ap:mods}
We make the following modifications to the feeder models to improve the realism of our study. 
 \begin{itemize}
  	\item Capacitor banks: By default, capacitor banks with more than one phase are configured so that the voltage-sensing phase is the only phase that is controlled. We modify these types of banks so that all phases are controllable. 
 	\item Power factors: The disaggregation method \cite{pop_script} sets houses' zip-loads to a power factor of 1.0. We adjust the power factors using Table A.2 of \cite{cohen_effects_2016} to better represent common loads in residential and commercial buildings.
 	\item Small static loads: In the original feeder models, there are a number of small static loads. By default, method \cite{pop_script} converts these loads to ``street lights'', which are only on in the evening. We revert these back to constant loads so that the loads are on during the peak hour.
 	\item Transformer sizing: The disaggregation method \cite{pop_script} can result in some transformers being undersized. We increase the size of transformers whose average loading is both greater than the transformer's original rating and greater than its original planning load by selecting the next largest transformer.
 \end{itemize}

\section{Cycling Behavior of Controlled ACs} \label{appendix:duty}

We show that, in a heterogeneous population of ACs under probabilistic dispatch in which all ACs receive the same $u$, some ACs will be more likely to be dispatch-switched off than on, and others will be more likely to be dispatch-switched on than off. We then show that the natural duty cycle of an AC can be used to predict which behavior will be most likely.

When an AC is on, the probability of it dispatch-switching off is equal to one minus the probability of it not dispatch-switching off during the same time. Thus, the probability of an AC dispatch-switching off at least once during the time it would take to complete the on part of its natural cycle is 
\begin{equation} \label{eq:Poff}
P_\text{S,OFF} = 1-\prod_{k=1}^{T_\text{ON}} (1-u_\text{OFF}(k)),\\
\end{equation}
and the probability of externally switching on at least once during the off part of its natural cycle is
\begin{equation} \label{eq:Pon}
P_\text{S,ON} = 1-\prod_{k=1}^{T_\text{OFF}} (1-u_\text{ON}(k)).\\
\end{equation}
Here $T_\text{ON}$ and $T_\text{OFF}$ are the number of time steps the AC is on and off during a natural cycle, and $u_\text{OFF}$ and $u_\text{ON}$ are the probabilistic dispatch commands for the off and on directions. 

If an AC has $P_\text{S,ON} > P_\text{S,OFF}$, then it is more likely to dispatch-switch on than off; if $P_\text{S,ON} < P_\text{S,OFF}$, the AC is more likely to dispatch-switch off than on. To determine an AC's behavior a priori, we estimate $P_\text{S,ON}$ and $P_\text{S,OFF}$ with a few simplifying assumptions. Let $D$ be the average of the ACs' natural duty cycles and $N$ be the total number of ACs. We assume that the regulation signal's amplitude is relatively small such that the ratio of the number of ACs on to number of ACs off remains approximately equal to $D$. Thus, if we switch a fixed number of ACs $N_\text{S}$ in each time step, we can approximate the switching commands as constant values: $u_\text{ON} \approx N_\text{S}/((1-D)N)$ and $u_\text{OFF} \approx N_\text{S}/(DN)$. Finally, let $d$ be the natural duty cycle of the given AC and $T_\text{P}$ be the number of time steps in its period, then $T_\text{ON} = dT_\text{P}$ and $T_\text{OFF}= (1-d)T_\text{P}$. With these substitutions, the probabilities of dispatch-switching during a natural cycle are 
\begin{equation} \label{eq:Poff2}
P_\text{S,OFF} \approx 1-\bigg(1-\frac{N_\text{S}}{DN}\bigg)^{dT_\text{P}} \text{ and}
\end{equation}
\begin{equation} \label{eq:Pon2}
P_\text{S,ON} \approx 1-\bigg(1-\frac{N_\text{S}}{(1-D)N}\bigg)^{(1-d)T_\text{P}}.
\end{equation}
After manipulating \eqref{eq:Poff2} and \eqref{eq:Pon2}, we find that an AC is more likely to dispatch-switch on than off (i.e., $P_\text{S,ON}>P_\text{S,OFF}$) if
\begin{equation}\label{eq:dutycycle}
\frac{d}{(1-d)} <  \frac{\log\Big(1-\frac{N_\text{S}}{(1-D)N}\Big)}{\log\Big(1-\frac{N_\text{S}}{DN}\Big)}. 
\end{equation}
Note that if $D = 0.5$ then the right-hand side of \eqref{eq:dutycycle} is equal to one. In this case, an individual AC will be more likely to be dispatch-switch on than off if $d<D$, which is what we observe empirically. 
\end{appendices}

\section*{Acknowledgment}

The authors thank David Chassin, Frank Tuffner, and Jason Fuller for their help with GridLAB-D.

\bibliographystyle{IEEEtran}

\bibliography{sjc_references}

\end{document}